\documentclass[a4paper,twoside]{article}
\usepackage{natbib}
\usepackage{amsmath,amsfonts,amssymb,amscd,graphics}
\usepackage{theorem}
\usepackage[english]{babel}

\theoremstyle{plain}
\newtheorem{theo}{Theorem}
\newtheorem{definition}[theo]{Definition}

\newtheorem{prop}[theo]{Proposition}
\newtheorem{lemme}[theo]{Lemma}
\newtheorem{cor}[theo]{Corollary}

{\theorembodyfont{\upshape}
}

{\theorembodyfont{\upshape}
}

{\theorembodyfont{\upshape}
\newtheorem{rem}[theo]{Remark}}

\def\demo{\noindent \textsc{Proof} \\}
\def\findem{\hfill$\square$\vskip 13pt}

\def\N{\mathbb{N}}
\def\R{\mathbb{R}}

\def\E{\mathbb{E}}
\def\F{\mathcal{F}}
\def\G{\mathcal{G}}
\def\B{\mathcal{B}}
\def\P{\mathbb{P}}
\def\cvp{\xrightarrow{\P}}
\def\cvl{\xrightarrow{\mathcal{L}}}
\def\cvps{\xrightarrow{a.s.}}
\def\cvBC{\xrightarrow{BC}}
\def\cvL1{\xrightarrow{L^1}}

\title{Convergence of values in optimal stopping \\
and convergence of optimal stopping times}

\author{Fran\c cois COQUET$^*$ and Sandrine TOLDO$^{**}$
\medskip\\
$^*${\it CREST-ENSAI, Campus de Ker Lann, 35170 Bruz, France}\\
{\it and LMAH, Universit\'e du Havre, France. E-mail:} {\rm fcoquet@ensai.fr} \\
$^{**}${\it IRMAR, Antenne de Bretagne de l'ENS Cachan, Campus de Ker Lann,}\\
{\it 35170 Bruz, France.  E-mail:} {\rm sandrine.toldo@bretagne.ens-cachan.fr}}

\date{}

\begin{document}

\maketitle

\noindent \textbf{Abstract:} Under the hypothesis of convergence in
probability of a sequence of c\`adl\`ag processes $(X^n)_n$ to a
c\`adl\`ag process $X$, we are interested in the convergence of
corresponding values in optimal stopping and also in the convergence
of optimal
stopping times. We give results under hypothesis of inclusion of filtrations or convergence of filtrations. \\


\section{Introduction}
\label{intro}

Let us consider a c\`adl\`ag process $X$. Denote by $\F^X$ its natural filtration and by $\F$ the right-continuous associated
filtration ($\forall t, \F_t = \F^X_{t^+}$). Denote also by $\mathcal{T}_T$ the set of $\F-$stopping times bounded by $T$. \\
\indent Let $\gamma: [0, T ] \times \R \to \R$ a bounded continuous
function. We define the value in optimal stopping of horizon $T$ for
the process $X$ by:
$$\Gamma(T)=\underset{\tau \in \mathcal{T}_T}{\sup}\E[\gamma(\tau,X_\tau)].$$
We call a stopping time $\tau$ {\it optimal} whenever $\E[\gamma(\tau,X_\tau)]=\Gamma(T)$.

\begin{rem}
\label{GammaX}
As it is noticed in \citet{cvreduites}, the value of $\Gamma(T)$ only depends on the law of $X$.
\end{rem}

Throughout this paper, we will deal with the problem of stability of
values in optimal stopping, and of optimal stopping times, under
approximations of the process $X$. To be more precise, let us
consider a sequence $(X^n)_n$ of processes which converges in
probability to a limit process $X$. For all $n$, we denote by $\F^n$
the natural filtration of $X^n$ and by $\mathcal{T}^n_T$ the set of
$\F^n-$stopping times bounded by $T$. Then, we define the values in
optimal stopping $\Gamma_n(T)$ by $\Gamma_n(T)=\underset{\tau \in
\mathcal{T}^n_T}{\sup}\E[\gamma(\tau,X^n_\tau)].$ The main aims of
this paper are first to give conditions under which
$(\Gamma_n(T))_n$ converges to $\Gamma(T)$, and second, when it is
possible to find a sequence $(\tau_n)$ of optimal stopping
times w.r.t. the $X_n$'s, to give further conditions under which the sequence $(\tau_n)$ converges to an optimal stopping time w.r.t $X$. \\

\indent In his unpublished manuscript \citep{preprintAldous}, Aldous proved that if $X$ is quasi-left continuous and if extended
convergence (in law) of $((X^n,\F^n))_n$ to $(X,\F)$ holds, then $(\Gamma_n(T))_n$ converges to $\Gamma(T)$. In their paper
\citep{cvreduites}, Lamberton and Pag\`es obtained the same result under the hypothesis of weak extended convergence of
$((X^n,\F^n))_n$ to $(X,\F)$, quasi-left continuity of the $X^n$'s and Aldous' criterion of tightness for $(X^n)_n$. \\
\indent Another way to study this problem is to consider the Snell envelopes associated to the processes. We recall that the Snell
envelope of a process $Y$ is the smallest supermartingale larger than $Y$ (see e.g. \citep{Karoui}). The value in optimal stopping can be
written as the value at 0 of a Snell envelope, as it is used for example in \citep{MP}, where a result of convergence of Snell envelopes
for the Meyer-Zheng topology is proved. \\

\indent Section \ref{cvvaleurs} is devoted to convergence of values in optimal stopping. The main difficulty is to prove that $\Gamma(T)
\geqslant \limsup \Gamma_n(T)$ and both papers \citep{preprintAldous} and \citep{cvreduites} need weak extended convergence to prove it.
We prove that this inequality actually holds whenever filtrations $\F^n$ are included into the limiting filtration $\F$, or when convergence
of filtrations holds. \\
\indent The main idea in our proof of the inequality $\Gamma(T) \geqslant \limsup \Gamma_n(T)$ is the following. We build a sequence
$(\tau^n)$ of $\F^n-$stopping times bounded by $T$. Then, we extract a convergent subsequence of $(\tau^n)$ to a random variable
$\tau$ and, at the same time, we compare $\E[\gamma(\tau,X_\tau)]$ and $\Gamma(T)$. This is carried out through two methods.\\
\indent First, we enlarge the space of stopping times, by considering the ran\-do\-mi\-zed stopping times and the topology
introduced in \citep{BC}. Baxter and Chacon have shown that the space of randomized stopping times with respect to a
right-continuous filtration with the associated topology is compact. We use this method in subsection \ref{limsup_inclusion} when holds
the hypothesis of inclusion of filtrations $\F^n \subset \F$ (which means that $\forall t \in [0,T], \F^n_t \subset \F_t$). We point out
that this assumption is simpler and easier to check than the extended convergence used in \citep{preprintAldous} and \citep{cvreduites}
or our own alternate hypothesis of convergence of filtrations. \\
\indent However, when inclusion of filtrations does not hold, we follow an idea already used, in a slightly different way, in
\citep{preprintAldous} and in \citep{cvreduites}, that is to enlarge the filtration $\F$ associated to the limiting process $X$. In
subsection \ref{limsup_cvfiltrations}, we enlarge (as little as possible) the limiting filtration so that the limit $\tau^*$ of a
convergent subsequence of the randomized $(\F^n)$ stopping times associated to the $(\tau^n)_n$ is a randomized stopping time for
this enlarged filtration and we assume that convergence of filtrations (but not necessarily extended convergence) holds. In doing so, we do
not need to introduce the prediction process which Aldous needed to define extended convergence. We also point out that convergence of
processes joined to convergence of filtrations does not always imply extended convergence (see \citep{Memin_2003} for a counter example). So
the result given in this subsection is somewhat different from those of \citep{preprintAldous} and \citep{cvreduites}. \\

When convergence of values in optimal stopping holds, it is natural 
to wonder wether the associated optimal stopping times (when 
existing) do converge. Here again, the main problem is that, in 
general, the limit of a sequence of stopping times is not a stopping 
time. It may happen that the limit in law of a sequence of 
$\F-$stopping times is not the law of a $\F-$stopping time (see the 
example in \citep{BC}). In section \ref{cvtaopt}, we shall give 
conditions, including again convergence of filtrations, under which 
the limit in probability of a sequence of $(\F^n)-$stopping times is 
a $\F-$stopping time (and not only a stopping time for a larger 
filtration). This caracterization will allow us to deduce a result 
of convergence of optimal stopping times when the limit process $X$ 
has independent increments. \\

\indent Finally, in section \ref{appl}, we give applications of the previous results to discretizations and also to financial models. \\

In what follows, we are given a probability space $(\Omega, \mathcal{A}, \P)$. We fix a positive real $T$.
Unless otherwise specified, every $\sigma$-field is supposed to be included in $\mathcal{A}$, every process will be indexed by $[0,T]$ and
taking values in $\R$ and every filtration will be indexed by $[0,T]$. $\mathbb{D}=\mathbb{D}([0,T])$ denotes the space of c\`adl\`ag
functions from $[0,T]$ to $\R$. We endow $\mathbb{D}$ with the Skorokhod topology. \\
\indent For technical background about Skorokhod topology, the reader may refer to \citep{Bill} or \citep{JS}.\\


\section{Convergence of values in optimal stopping}\label{cvvaleurs}
\subsection{Statement of the results}

The notion of convergence of filtrations has been defined in
\citep{Hoover} and, in a slightly different way, in
\citep{cvfiltration}. Here, we use the definition taken from the
latter paper:

\begin{definition}
We say that $(\F^n)$ converges weakly to $\F$ if for every $A \in \F_T$, $(\E[1_A | \F^n_.])_n$
converges in probability to $\E[1_A | \F_.]$ for the Skorokhod topology. We denote $\F^n \xrightarrow{w} \F$.
\end{definition}

Aldous' Criterion for tightness,which has been introduced in the
papers \citep{Aldous78} and \citep{Aldous89}, is a standard tool for
functional limit theorems when the limit is quasi-left continuous.
It happens to be at the heart of the following Theorem, whose proof
is the main purpose of this section.

\begin{theo}
\label{cvGamma} Let us consider a c\`adl\`ag process X and a
sequence $(X^n)_n$ of c\`adl\`ag processes. Let $\F$ be the
right-continuous filtration associated to the natural filtration of
X and $(\F^n)_n$ the natural filtrations of the processes $(X^n)_n$.
We assume that $X^n \cvp X$, that Aldous' Criterion for tightness is
filled, i.e.
\begin{equation}
\label{hypA}
\forall \varepsilon > 0, \underset{\delta \downarrow 0}{\lim} \ \underset{n \to +\infty}{\limsup} \
        \underset{\sigma, \nu \in \mathcal{T}_T^n, \sigma \leqslant \nu \leqslant \sigma+\delta}{\sup} \
                \P[| X^n_{\sigma} - X^n_{\nu} | \geqslant \varepsilon] = 0,
\end{equation}
and that one of the following assertions holds: \\
- for all n, $\F^n \subset \F$, \\
- $\F^n \xrightarrow{w} \F$. \\
Then, $\Gamma_n(T) \xrightarrow[n \to \infty]{} \Gamma(T)$.
\end{theo}

The proof of Theorem \ref{cvGamma} will be carried out through two steps in next subsections: \\
- Step 1: show that $\Gamma(T) \leqslant \liminf \Gamma_n(T)$ in subsection \ref{liminf}, \\
- Step 2: show that $\Gamma(T) \geqslant \limsup \Gamma_n(T)$ in subsections \ref{limsup_inclusion} and \ref{limsup_cvfiltrations}. \\

Let us give at once an extension of Theorem \ref{cvGamma} which will prove useful for the application to finance in Section \ref{appl}.
\begin{cor}
\label{cvGamma^n}
Let $(\gamma^n)_n$ be a sequence of continuous bounded functions on $[0,T] \times \R$ which uniformly converges to a continuous bounded
function $\gamma$. Let $X$ be a c\`adl\`ag process and $(X^n)_n$ a sequence of c\`adl\`ag processes. Let $\F$ be the right-continuous
filtration of the process $X$ and $\F^n$ the natural filtration of $X^n$.
We suppose that $X^n \cvp X$, that Aldous' Criterion for tightness (\ref{hypA}) is filled and that one of the following assertions holds: \\
- for every $n$, $\F^n \subset \F$, \\
- $\F^n \xrightarrow{w} \F$. \\
We consider the values in optimal stopping defined by:
$$\Gamma(T)=\underset{\tau \in \mathcal{T}_T}{\sup}\E[\gamma(\tau,X_\tau)] \text{~~and~~}
    \Gamma_n(T)=\underset{\tau^n \in \mathcal{T}^n_T}{\sup}\E[\gamma^n(\tau^n,X^n_{\tau^n})].$$
Then $\Gamma_n(T) \xrightarrow[n \to \infty]{} \Gamma(T)$.
\end{cor}

\demo
According to Theorem \ref{cvGamma}, we have
$\underset{\tau^n \in \mathcal{T}^n_T}{\sup}\E[\gamma(\tau^n,X^n_{\tau^n})]
    \to \underset{\tau \in \mathcal{T}_T}{\sup}\E[\gamma(\tau,X_\tau)].$
To conclude, it suffices to prove that
$\underset{\tau^n \in \mathcal{T}^n_T}{\sup}\E[\gamma^n(\tau^n,X^n_{\tau^n})]
    - \underset{\tau^n \in \mathcal{T}^n_T}{\sup}\E[\gamma(\tau^n,X^n_{\tau^n})]
        \to 0.$
But,
\begin{eqnarray*}
\underset{\tau^n \in \mathcal{T}^n_T}{\sup}\E[\gamma^n(\tau^n,X^n_{\tau^n})]
    - \underset{\tau^n \in \mathcal{T}^n_T}{\sup}\E[\gamma(\tau^n,X^n_{\tau^n})]
& \leqslant & \underset{\tau^n \in \mathcal{T}^n_T}{\sup} |\E[\gamma^n(\tau^n,X^n_{\tau^n})] - \E[\gamma(\tau^n,X^n_{\tau^n})]| \\
& \leqslant & \underset{\tau^n \in \mathcal{T}^n_T}{\sup} \E[|\gamma^n(\tau^n,X^n_{\tau^n}) - \gamma(\tau^n,X^n_{\tau^n})|] \\
& \leqslant & \sup_{t,x} |\gamma^n(t,x) - \gamma(t,x)| \\
& \to & 0 \text{~~using the uniform convergence.}
\end{eqnarray*}
Corollary \ref{cvGamma^n} is proved.
\findem


\subsection{Proof of the inequality $\Gamma(T) \leqslant \liminf \Gamma_n(T)$}
\label{liminf}

In this section, we give a lower semi-continuity result. The hypotheses are not the weakest possible, but will be sufficient to
prove Theorem \ref{cvGamma}.

\begin{theo}
\label{thliminf}
Let us consider a c\`adl\`ag process $X$, the right-continuous filtration $\F$ associated to the natural filtration of $X$, a sequence of
c\`adl\`ag processes $(X^n)_n$ and their natural filtrations $(\F^n)_n$. We suppose that \mbox{$X^n \cvp X$.}
Then $\Gamma(T) \leqslant \liminf \Gamma_n(T)$.
\end{theo}

\demo
We only give here the sketch of the proof, which is not very different from those in \citep{cvreduites} and \citep{preprintAldous}. \\

To begin with, we can prove that, if $\tau$ is a $\F^X-$stopping
time bounded by $T$ and taking values in a discrete set $\{t_i\}_{i
\in I}$ such that $\P[\Delta X_{t_i} \not= 0]=0$, $\forall i$, and
if we define $\tau^n$ by \mbox{$\tau^n(\omega)=\min \{t_i: i \in\{j
: \E[1_{A_j}|\F^n_{t_j}](\omega) > 1/2\}\}$}, $\forall \omega$,
where \mbox{$A_i = \{\tau=t_i\}$}, then, $(\tau^n)$ is a sequence of
$(\mathcal{T}^n_T)$ such that $(\tau^n,X^n_{\tau^n}) \cvp
(\tau,X_\tau)$.

Let us then consider a finite subdivision $\pi$ of $[0,T]$ such that
no fixed time of discontinuity of $X$ belongs to $\pi$. We denote by
$\mathcal{T}^{\pi}_T$ the set of $\F$ stopping times taking values
in $\pi$, and we define:
$$\Gamma^{\pi}(T)=\underset{\tau \in \mathcal{T}^{\pi}_T}{\sup}\E[\gamma(\tau,X_\tau)].$$
Applying the previous result to stopping times belonging to $\mathcal{T}^{\pi}_T$ shows that
$\Gamma^{\pi}(T) \leqslant \liminf  \Gamma_n(T).$

At last, using an increasing sequence $(\pi^k)_k$ of subdivisions such that $|\pi^k| \xrightarrow[k \to +\infty]{} 0$ and such that
$\P[\Delta X_s \not= 0]=0$ $\forall s \in \pi_k$, standard computations prove that $\Gamma^{\pi^k}(T) \xrightarrow[k \to +\infty]{} \Gamma(T)$,
and Theorem \ref{thliminf} follows.
\findem


\subsection{Proof of the inequality $\Gamma(T) \geqslant \limsup \Gamma_n(T)$ if for every $n$, $\F^n \subset \F$}
\label{limsup_inclusion}

\subsubsection{Randomized stopping times }
\label{tar}

The notion of randomized stopping times has been introduced in \citep{BC} and this notion has been used in \citep{Meyer}
 under the french name "temps d'arr\^et flous". \\

We are given a filtration $\F$. Let us denote by $\mathcal{B}$ the
Borel $\sigma$-field on $[0,1]$. Then, we define the filtration
$\mathcal{G}$ on $\Omega \times [0,1]$ such that $\forall t$,
$\mathcal{G}_t=\F_t \times \mathcal{B}$. A map \mbox{$\tau: \Omega
\times [0,1] \to [0,+\infty]$} is called a randomized $\F-$stopping
time if $\tau$ is a $\mathcal{G}-$stopping time. We denote by
$\mathcal{T}^*$ the set of randomized stopping times and by
$\mathcal{T}^*_T$ the set of randomized stopping times bounded by
$T$. $\mathcal{T}$ is included in $\mathcal{T}^*$ and the
application $\tau \mapsto \tau^*$, where
$\tau^*(\omega,t)=\tau(\omega)$ for every $\omega$ and every $t$,
maps $\mathcal{T}$ into $\mathcal{T}^*$. In the same way, 
$\mathcal{T}_T$ is included in $\mathcal{T}^*_T$. \\ \indent On the 
space $\Omega \times [0,1]$, we build the probability measure $\P 
\otimes \mu$ where $\mu$ is Lebesgue's measure on $[0,1]$. In their
paper \citep{BC}, Baxter and Chacon define the convergence of 
randomized stopping times by the following: $$\tau^{*,n} 
\xrightarrow{BC} \tau^* \text{~~iff~~} \forall f \in \mathcal{C}_b
([0,\infty]), \forall Y \in L^1(\Omega, \F, \P), \E[Yf(\tau^{*,n})]
\to \E[Yf(\tau^*)],$$ where $\mathcal{C}_b ([0,\infty])$ is the set
of bounded continuous functions on $[0,\infty]$.\\
\indent This kind of convergence
is a particular case of "stable convergence" as introduced in
\citep{Renyi} and studied in \citep{JacodMemin}. \\

\indent The main point for us here is, as it is shown in \citep[Theorem 1.5]{BC}, that the set of randomized stopping times for a
right-continuous filtration is compact for Baxter and Chacon's topology (which is not true for the set of ordinary stopping times). \\

The following Proposition will be the main argument in the proof of Theorem \ref{thlimsup1} below.

\begin{prop}
\label{sssuitetar}
Let us consider a sequence of filtrations $(\F^n)$ and a right-continuous filtration $\F$ such that $\forall n$, $\F^n \subset \F$.
Let $(\tau^n)_n$ be a sequence of $(\mathcal{T}_T^{n})_n$. Then, there exists a randomized $\F-$stopping time $\tau^*$ and a subsequence
$(\tau^{\varphi(n)})_n$ such that $\tau^{*,\varphi(n)} \xrightarrow{BC} \tau^*$ where for every $n$, $\tau^{*,n}(\omega,t)=\tau^n(\omega)$
$\forall \omega$, $\forall t$.
\end{prop}

\demo
For every $n$, $\F^n \subset \F$, $(\tau^n)_n$ is a sequence of $\F-$stopping times so, by definition, $(\tau^{*,n} )_n$ is a sequence of
randomized $\F-$stopping times. According to \citep[Theorem 1.5]{BC}, we can find a randomized $\F-$stopping time $\tau^*$ and a
subsequence $(\tau^{\varphi(n)})_n$ such that $\tau^{*,\varphi(n)} \xrightarrow{BC} \tau^*$.
\findem

\indent Now, we define $X_{\tau^*}$ by
$X_{\tau^*}(\omega,v)=X_{\tau^*(\omega,v)}(\omega)$, for every
$(\omega,v) \in \Omega \times [0,1]$. Then, we can prove the
following Lemma:
\begin{lemme}
\label{GammaGamma*} Let $\Gamma^*(T)=\underset{\tau^* \in
\mathcal{T}^*_T}{\sup}\E[\gamma(\tau^*,X_{\tau^*})]$. Then
$\Gamma^*(T)=\Gamma(T)$.
\end{lemme}

\demo
- $\mathcal{T}_T$ is included into $\mathcal{T}^*_T$, hence $\Gamma(T) \leqslant \Gamma^*(T)$. \\
- Let $\tau^* \in \mathcal{T}^*_T$. We consider, for every $v$, $\tau_v(\omega)=\tau^*(\omega,v), \forall \omega$. \\
For every $v \in [0,1]$, for every $t \in [0,T]$,
$$\{\omega: \tau_v(\omega) \leqslant t\} \times \{v\} = \{ (\omega, x): \tau^*(\omega,x) \leqslant t\} \cap (\Omega \times \{v\}).$$
But, $\{ (\omega, x): \tau^*(\omega,x) \leqslant t\} \in \F_t \times
\mathcal{B}$ because $\tau^*$ is a randomized $\F-$stopping time and
$\Omega \times \{v\} \in \F_t \times \mathcal{B}$. So, $\{\omega:
\tau_v(\omega) \leqslant t\} \times \{v\} \in \F_t \times
\mathcal{B}$. Consequently,
$$\{\omega: \tau_v(\omega) \leqslant t\} \in \F_t.$$
Hence, for every $v$, $\tau_v$ is a $\F-$stopping time bounded by
$T$. We have:
\begin{eqnarray*}
\E[\gamma(\tau^*,X_{\tau^*})] & = & \int_\Omega \int_0^1 \gamma(\tau^*(\omega,v),X_{\tau^*(\omega,v)}(\omega))d\P(\omega)dv \\
& = & \int_0^1 \left( \int_\Omega \gamma(\tau^*(\omega,v),X_{\tau_v(\omega)}(\omega))d\P(\omega) \right) dv \\
& = & \int_0^1 \E[\gamma(\tau_v,X_{\tau_v})] dv \\
& \leqslant & \Gamma(T) \text{~~because, for every $v$, $\tau_v \in \mathcal{T}_T$.}
\end{eqnarray*}
Taking the $\sup$ over $\tau^*$ in $\mathcal{T}^*_T$, we get $\Gamma^*(T) \leqslant \Gamma(T)$. \\
Lemma \ref{GammaGamma*} is proved.
\findem

The following proposition will also be useful:

\begin{prop}
\label{cvlXntar}
Let us consider a sequence $(X^n)_n$ of c\`adl\`ag adapted processes that converges in law to a c\`adl\`ag process X. Let $(\tau^n)_n$ be
a sequence of stopping times such that the associated sequence $(\tau^{*,n})_n$ of randomized stopping times
($\tau^{*,n}(\omega,t)=\tau^n(\omega)$ $\forall \omega$, $\forall t$) converges in law to a random variable $V$. We suppose that
$(\tau^{*,n},X^n) \cvl (V,X)$ and that Aldous' Criterion (\ref{hypA}) is filled. Then $(\tau^{*,n},X^n_{\tau^{*,n}}) \cvl (V,X_V)$.
\end{prop}

\begin{rem}
As the proof of Proposition \ref{cvlXntar} follows the lines of the proof of \citep[Corollary 16.23]{preprintAldous}, we skip it here.
However we point out that, in this proposition, Aldous' Criterion is filled for genuine -not randomized- stopping times.
\end{rem}


\begin{prop}
\label{cvcouple}
Let us consider a sequence $(X^n)_n$ of c\`adl\`ag adated processes converging in probability to a c\`adl\`ag process X.
Let $(\tau^{*,n})_n$ be a sequence of ran\-do\-mi\-zed stopping times converging to the ran\-do\-mi\-zed stopping time $\tau$ under Baxter
and Chacon's topology. \\
Then $(X^n, \tau^{*,n}) \cvl (X, \tau^*)$.
\end{prop}

\demo
- As $(X^n)_n$ and $(\tau^{*,n})_n$ are tight, $((X^n, \tau^{*,n}))_n$ is tight for the product topology. \\
- We are now going to identify the limit throught finite-dimensional convergence. \\
Let $k \in \N$ and $t_1 < \ldots < t_k$ such that for every $i$, $\P[\Delta X_{t_i} \not= 0] = 0$.
We are going to show that $(X^n_{t_1}, \ldots , X^n_{t_k},\tau^{*,n}) \cvl (X_{t_1}, \ldots , X_{t_k}, \tau^*)$. \\
Let $f: \R^k \to \R$ and $g: \R \to \R$ be bounded continuous
functions.
\begin{eqnarray*}
\lefteqn{ |\E[f(X^n_{t_1}, \ldots , X^n_{t_k})g(\tau^{*,n})] - \E[f(X_{t_1}, \ldots , X_{t_k})g(\tau^*)]| } \\
& \leqslant & |\E[(f(X^n_{t_1}, \ldots , X^n_{t_k})-f(X_{t_1}, \ldots , X_{t_k}))g(\tau^{*,n})]| \\
& & \quad + |\E[f(X_{t_1}, \ldots , X_{t_k})g(\tau^{*,n})]-\E[f(X_{t_1}, \ldots , X_{t_k})g(\tau^*)]| \\
& \leqslant & \|g\|_{\infty} \E[|f(X^n_{t_1}, \ldots , X^n_{t_k})-f(X_{t_1}, \ldots , X_{t_k})|] \\
& & \quad + |\E[f(X_{t_1}, \ldots , X_{t_k})g(\tau^{*,n})]-\E[f(X_{t_1}, \ldots , X_{t_k})g(\tau^*)]| \\
\end{eqnarray*}
But, $X^n \cvp X$ and for every $i$, $\P[\Delta X_{t_i} \not= 0] = 0$ so
\mbox{$(X^n_{t_1}, \ldots , X^n_{t_k}) \cvp (X_{t_1}, \ldots , X_{t_k})$.}
Moreover, $f$ is bounded continuous, so
$$\E[|f(X^n_{t_1}, \ldots , X^n_{t_k})-f(X_{t_1}, \ldots , X_{t_k})|] \xrightarrow[n \to +\infty]{} 0.$$
On the other hand, by definition of Baxter and Chacon's convergence,
$$\E[f(X_{t_1}, \ldots , X_{t_k})g(\tau^{*,n})]-\E[f(X_{t_1}, \ldots , X_{t_k})g(\tau^*)] \xrightarrow[n \to +\infty]{} 0.$$
Then,
$$\E[f(X^n_{t_1}, \ldots , X^n_{t_k})g(\tau^{*,n})] - \E[f(X_{t_1}, \ldots , X_{t_k})g(\tau^*)] \xrightarrow[n \to +\infty]{} 0.$$
\indent Using a density argument, we can expand the previous result
to continuous and bounded fonctions from $\R^{k+1}$ to $\R$. More
precisely, for every $\varphi: \R^{k+1} \to \R$ continuous and
bounded, we have:
$$\E[\varphi(X^n_{t_1}, \ldots , X^n_{t_k},\tau^{*,n})] \xrightarrow[n \to +\infty]{} \E[\varphi(X_{t_1}, \ldots , X_{t_k},\tau^*)].$$
It follows that $(X^n_{t_1}, \ldots , X^n_{t_k},\tau^{*,n}) \cvl (X_{t_1}, \ldots , X_{t_k}, \tau^*)$. At last, the tightness of the
sequence $((X^n, \tau^{*,n}))_n$ and the finite-dimensional convergence on a dense set to $(X, \tau^*)$ imply
$(X^n, \tau^{*,n}) \cvl (X,\tau^*)$.
\findem

\subsubsection{Application to the proof of the inequality $\limsup \Gamma_n(T) \leqslant \Gamma(T)$}

We can now prove our first result about convergence of optimal
values.

\begin{theo}
\label{thlimsup1}
Let us consider a c\`adl\`ag process X, its right-continuous filtration $\F$, a sequence $(X^n)_n$ of
c\`adl\`ag processes and their natural filtrations $(\F^n)_n$. We suppose that $X^n \cvp X$, that Aldous' Criterion for tightness
(\ref{hypA}) is filled and that $\forall n$, $\F^n \subset \F$. Then $\limsup \Gamma_n(T) \leqslant \Gamma(T)$.
\end{theo}

\demo There exists a subsequence $(\Gamma_{\varphi(n)}(T))_n$ which 
converges to $\limsup \Gamma_n(T)$. \\ Let us fix $\varepsilon > 0$. 
We can find a sequence $(\tau^{\varphi(n)})_n$ of 
$(\mathcal{T}^{\varphi(n)}_T)_n$ such that $$\forall n,\ 
\E[\gamma(\tau^{\varphi(n)}, X^{\varphi(n)}_{\tau^{\varphi(n)}})]
\geqslant \Gamma_{\varphi(n)}(T) - \varepsilon.$$ We consider the
sequence $(\tau^{*,\varphi(n)})_n$ of randomized stopping times 
associated to \mbox{$(\tau^{\varphi(n)})_n$:} for every $n$, 
$\tau^{*,\varphi(n)}(\omega,t)=\tau^{\varphi(n)}(\omega)$, $\forall 
\omega$, $\forall t$. $\F^{\varphi(n)} \subset \F$ and
$(\tau^{\varphi(n)})$ is a sequence of $\F^{\varphi(n)}-$stopping 
times bounded by $T$, so using Proposition \ref{sssuitetar}, there 
exists a randomized $\F-$stopping time $\tau^*$ and a subsequence 
$(\tau^{\varphi \circ \psi(n)})$ such that $\tau^{*,\varphi \circ 
\psi(n)} \xrightarrow{BC} \tau^*$. Moreover $X^{\varphi \circ
\psi(n)} \cvp X$, so by Proposition \ref{cvcouple}, $(X^{\varphi
\circ \psi(n)}, \tau^{*,\varphi \circ \psi(n)}) \cvl (X,\tau^*).$
Then, using Proposition \ref{cvlXntar}, we have: $(\tau^{*,\varphi 
\circ \psi (n)},X^{\varphi \circ \psi (n)}_{\tau^{*,\varphi \circ
\psi (n)}}) \cvl (\tau^*,X_{\tau^*}).$ Since $\gamma$ is continuous
and bounded, we deduce: $$\E[\gamma(\tau^{*,\varphi \circ \psi (n)}, 
X^{\varphi \circ \psi(n)}_{\tau^{*,\varphi \circ \psi(n)}})] \to
\E[\gamma(\tau^*,X_{\tau^*})].$$ But, $\E[\gamma(\tau^{*,\varphi
\circ \psi(n)}, X^{\varphi \circ \psi(n)}_{\tau^{*,\varphi \circ
\psi(n)}})]
    = \E[\gamma(\tau^{\varphi \circ \psi(n)},X^{\varphi \circ \psi(n)}_{\tau^{\varphi \circ \psi(n)}})]$
by definition of $(\tau^{*,n})$, and by construction of $\varphi$,
$\E[\gamma(\tau^{\varphi \circ \psi(n)},X^{\varphi \circ \psi(n)}_{\tau^{\varphi \circ \psi(n)}})]
    \geqslant \Gamma_{\varphi \circ \psi(n)}(T) - \varepsilon$.
So,
$$\E[\gamma(\tau^*,X_{\tau^*})] \geqslant \lim\Gamma_{\varphi \circ \psi(n)}(T) - \varepsilon=\limsup\Gamma_n(T)-\varepsilon.$$

We hence have proved that for any $\varepsilon>0$ we can find a randomized stopping time $\tau^*$ such that
$\E[\gamma(\tau^*,X_{\tau^*})] \geqslant \limsup\Gamma_n(T)-\varepsilon.$

As by definition $\E[\gamma(\tau^*,X_{\tau^*})] \leqslant \Gamma^*(T)$ and $\varepsilon$ is arbitrary, it follows that
$\Gamma^*(T) \geqslant \limsup \Gamma_n(T)$.

At last, recall that $\Gamma^*(T)=\Gamma(T)$ by Lemma \ref{GammaGamma*} to conclude that $\Gamma(T) \geqslant \limsup \Gamma_n(T)$.
\findem

\begin{rem}
We were able to prove the previous theorem, because we knew something about the nature of the limit of the subsequence of stopping
times thanks to Proposition \ref{sssuitetar}. If we remove the hypothesis of inclusion of filtrations $\F^n \subset \F, \forall n$, the
limit of the subsequence needs no longer be a randomized $\F-$stopping time, and we cannot always compare $\E[\gamma(\tau^*,X_{\tau^*})]$
to $\Gamma^*(T)$.

However, the result of Theorem \ref{thlimsup1} remains true under other settings, as we shall prove in next subsection.
\end{rem}


\subsection{Proof of the inequality $\limsup \Gamma_n(T) \leqslant \Gamma(T)$ if $\F^n \xrightarrow{w} \F$}
\label{limsup_cvfiltrations}

\begin{theo}
\label{thlimsup2}
Let us consider a sequence of c\`adl\`ag processes $(X^n)_n$, their na\-tu\-ral filtrations $(\F^n)_n$, a c\`adl\`ag process $X$
and its right-continuous natural filtration $\F$. We suppose that $X^n \cvp X$, that Aldous' Criterion for tightness (\ref{hypA}) is filled
and that $\F^n \xrightarrow{w} \F$. Then $\limsup \Gamma_n(T) \leqslant \Gamma(T)$.
\end{theo}

\demo Our proof is more or less scheduled as the second part of the
proof in \citep[Theorem 17.2]{preprintAldous}. The main difference
is that we do not need extended convergence in our theorem: instead, we use convergence of filtrations. \\
We can find a subsequence $(\Gamma_{\varphi(n)}(T))_n$ converging to $\limsup \Gamma_n(T)$. \\
Let us take $\varepsilon > 0$. There exists a sequence $(\tau^{\varphi(n)})_n$ of $(\mathcal{T}^{\varphi(n)}_T)_n$ such that
$$\forall n, \E[\gamma(\tau^{\varphi(n)}, X^{\varphi(n)}_{\tau^{\varphi(n)}})] \geqslant \Gamma_{\varphi(n)}(T) - \varepsilon.$$
Let us consider the sequence $(\tau^{*,{\varphi(n)}})_n$ of associated randomized $\F^{\varphi(n)}-$stopping times like in \ref{tar}.
Taking the filtration $\mathcal{H}=(\bigvee_n \F^n) \vee \F$, $(\tau^{*,{\varphi(n)}})$ is a bounded sequence of randomized
$\mathcal{H}-$stopping times. Then, using \citep[Theorem 1.5]{BC}, we can find a further subsequence (still denoted  $\varphi$) and a
randomized $\mathcal{H}-$stopping time $\tau^*$ ($\tau^*$ is not a priori a randomized $\F-$stopping time) such that
$$\tau^{*,\varphi(n)} \cvBC \tau^*.$$
Using Proposition \ref{cvcouple}, we obtain $(X^{\varphi(n)}, \tau^{*,\varphi(n)}) \cvl (X,\tau^*)$.
Then, Proposition \ref{cvlXntar} gives the convergence
$(\tau^{\varphi(n)},X^{\varphi(n)}_{\tau^{\varphi(n)}}) \cvl (\tau^*,X_{\tau^*}).$ So,
$\E[\gamma(\tau^{\varphi(n)},X^{\varphi(n)}_{\tau^{\varphi(n)}})]
        \xrightarrow[n \to +\infty]{} \E[\gamma(\tau^*,X_{\tau^*})].$
On the other hand, $\E[\gamma(\tau^{\varphi(n)},X^{\varphi(n)}_{\tau^{\varphi(n)}})]
        \geqslant \Gamma_{\varphi(n)}(T) - \varepsilon$. So, letting $n$ go to infinity leads to

\begin{equation}\label{limsup-epsilon}
\E[\gamma(\tau^*,X_{\tau^*})] \geqslant \limsup \Gamma_n(T) - \varepsilon.
\end{equation}

Our next step will be to prove the following
\begin{lemme}\label{mingamma}
$$\E[\gamma(\tau^*,X_{\tau^*})]\leqslant\Gamma^*(T).$$
\end{lemme}
\demo
Let us consider the smaller right-continuous filtration $\G$ such that $X$ is $\G-$adapted and $\tau^*$ is a randomized $\G-$stopping time.
It is clear that $\F \subset \G$.
For every $t$, we have
$$\G_t \times \mathcal{B}=\bigcap_{s > t} \sigma(A \times B ,\{\tau^* \leqslant u\}, A \in \F_s, u \leqslant s, B\in\mathcal{B}).$$

We consider the set $\tilde{\mathcal{T}}_T$ of randomized $\G-$stopping times bounded by $T$ and we define
$\tilde{\Gamma}(T)=\underset{\tilde{\tau} \in \tilde{\mathcal{T}}_T}{\sup}\E[\gamma(\tilde{\tau},X_{\tilde{\tau}})]$. \\

By definition of $\G$, $\tau^* \in \tilde{\mathcal{T}}_T$ so
\begin{equation}\label{tildestar} \E[\gamma(\tau^*,X_{\tau^*})] \leqslant \tilde{\Gamma}(T).
\end{equation}

In order to prove Lemma \ref{mingamma}, we will use the following
Lemma, which is an adaptation of \citep[Proposition 3.5]{cvreduites}
to our enlargement of filtration:
\begin{lemme}
\label{indepcond}
If $\G_t \times \mathcal{B}$ and $\F_T \times \mathcal{B}$ are conditionally independent given $\F_t \times \mathcal{B}$
for every $t \in [0,T]$, then $\tilde{\Gamma}(T)=\Gamma^*(T)$.
\end{lemme}

\demo
The proof is the same as the proof of \citep[Proposition 3.5]{cvreduites} with
$(\F_t \times \B)_{t \in [0,T]}$ and $(\G_t \times \B)_{t \in [0,T]}$ instead of $\F^Y$ and $\F$ and with a process $X^*$ such that
for every $\omega$, for every $v \in [0,1]$, for every $t \in [0,T]$,
$X^*_t(\omega,v)=X_t(\omega)$ instead of the process $Y$.
\findem

Back to Lemma \ref{mingamma}, we have to prove the conditional
independence required in Lemma \ref{indepcond} which, according to
\citep[Theorem 3]{BY}, is equivalent to the following assumption:
\begin{equation}
\label{indcond}
\forall t \in [0,T], \forall Z \in L^1(\F_T \times \mathcal{B}), \E[Z|\F_t \times \mathcal{B}]=\E[Z|\G_t \times \mathcal{B}].
\end{equation}

The main part of what is left in this subsection is devoted to show
that the assumptions of Theorem \ref{thlimsup2} do imply
(\ref{indcond}), therefore fulfilling the assumptions needed to make
Lemma \ref{indepcond} work. Note that in order to prove
(\ref{indcond}) \citep{preprintAldous} and in \citep{cvreduites} use
extended convergence, which needs not
hold under the hypothesis of Theorem \ref{thlimsup2} (see \citep{Memin_2003} for a counter-example). \\

\noindent Without loss of generality, we suppose from now on that $\tau^{*,n} \cvBC \tau^*$ instead of $\tau^{*,\varphi(n)} \cvBC \tau^*$. \\
Moreover, as $X^n \cvp X$ and Aldous' Criterion for tightness (\ref{hypA}) is filled, using the results of \citep{preprintAldous}, $X$ is
quasi-left continuous. \\

\noindent - As $\F \subset \G$, $\forall t\in[0,T]$, $\forall Z \in 
L^1(\F_T \times \mathcal{B})$, $\E_{\P \otimes \mu}[Z|\F_t \times
\mathcal{B}]$ is $\G_t \times \mathcal{B}-$measurable. \\ - We shall
show that $\forall t \in [0,T], \forall Z \in L^1(\F_T \times 
\mathcal{B}), \forall C \in \G_t \times \mathcal{B},$ $$\E_{\P
\otimes \mu}[\E_{\P \otimes \mu}[Z|\F_t \times
\mathcal{B}]1_C]=\E_{\P \otimes \mu}[Z1_C].$$

Let $t \in [0,T]$ and $\varepsilon > 0$ be fixed, and take $Z \in L^1(\F_T \times \mathcal{B})$.
By definition of $\G_t \times \mathcal{B}$, it suffices to prove that for every $A \in \F_t$, for every $s \leqslant t$ and for every
$B \in \mathcal{B}$,
\begin{eqnarray}
\label{eqZ}
\lefteqn{\int\int_{\Omega \times [0,1]} Z(\omega,v)1_A(\omega)1_{\{\tau^*(\omega,v) \leqslant s\}}1_B(v) d\P(\omega)dv }\\
& = &
\int\int_{\Omega \times [0,1]} \E_{\P \otimes \mu}[Z|\F_t \times \mathcal{B}](\omega,v)1_A(\omega)1_{\{\tau^*(\omega,v) \leqslant s\}}
                1_B(v) d\P(\omega)dv. \nonumber
\end{eqnarray}

We first prove that (\ref{eqZ}) holds for $Z=1_{A_1 \times A_2}$, $A_1 \in \F_T$, $A_2 \in \mathcal{B}$. \\
We can find $l \in \N$, $s_1 < \ldots < s_l$ and a continuous bounded function $f$ such that
\begin{equation}
\label{alpha1}
\E_{\P}[|1_{A_1}-f(X_{s_1}, \ldots, X_{s_l})|] \leqslant \varepsilon.
\end{equation}
Then
$$\int\int |1_{A_1 \times A_2}(\omega,v)-f(X_{s_1}(\omega), \ldots, X_{s_l}(\omega))1_{A_2}(v)|d\P(\omega)dv \leqslant \varepsilon.$$

Let us fix $A \in \F_t$.
We can find $k \in \N$, $t_1 < \ldots < t_k \leqslant t$ and $H:\R^k \to \R$ bounded continuous such
that
\begin{equation}
\label{alpha2}
\E_{\P}[|1_A-H(X_{t_1}, \ldots, X_{t_k})|] \leqslant \varepsilon.
\end{equation}
Let $u > t$ such that $\P[\Delta \E[f(X_{s_1}, \ldots, X_{s_l})|\F_u] \not= 0]=0$ and $\P[\tau^*=u]=0$. \\
Fix $s \leqslant t$.
We can find a bounded continuous function $G$ such that
\begin{equation}
\label{alpha3}
\E_{\P \otimes \mu}[|1_{\{\tau^* \leqslant s\}}-G(\tau^* \wedge u)|] \leqslant \varepsilon.
\end{equation}
$B \in \mathcal{B}$ and the set of continuous functions is dense
into $L^1(\mu)$, so there exists $g: \R \to \R$ bounded continuous
such that
\begin{equation}
\label{alpha45}
\int |1_{B}(v)-g(v)|dv \leqslant \varepsilon.
\end{equation}

\noindent We are going to show that
\begin{eqnarray*}
&& \int \int \E_{\P \otimes \mu}[f(X_{s_1}, \ldots, X_{s_l})1_{A_2}|\F_u \otimes \mathcal{B}](\omega,v)
                H(X_{t_1}(\omega), \ldots, X_{t_k}(\omega))\\
&& \quad \quad \quad G(\tau^*(\omega,v) \wedge u)g(v)d(\P \otimes \mu)(\omega,v)\\
& = & \int\int f(X_{s_1}(\omega), \ldots, X_{s_l}(\omega))1_{A_2}(v)H(X_{t_1}(\omega), \ldots, X_{t_k}(\omega))\\
&& \quad \quad \quad G(\tau^*(\omega,v) \wedge u)g(v)d(\P \otimes \mu)(\omega,v).
\end{eqnarray*}
$X^n \cvp X$ and $f$ is a bounded continuous function, so that
$$f(X^n_{s_1}, \ldots, X^n_{s_l}) \cvL1 f(X_{s_1}, \ldots, X_{s_l}).$$
Moreover, $\F^n \xrightarrow{w} \F$ so using \citep[Remark 2]{cvfiltration},
$$\E_{\P}[f(X^n_{s_1}, \ldots, X^n_{s_l})|\F^n_.] \cvp \E_{\P}[f(X_{s_1}, \ldots, X_{s_l})|\F_.].$$
Since $\P[\Delta \E[f(X_{s_1}, \ldots, X_{s_l})|\F_u] \not= 0]=0$, we have
$$\E_{\P}[f(X^n_{s_1}, \ldots, X^n_{s_l})|\F^n_u] \cvp \E_{\P}[f(X_{s_1}, \ldots, X_{s_l})|\F_u]$$
and since $f$ is bounded,
\begin{equation}
\label{mer1}
\E_{\P}[f(X^n_{s_1}, \ldots, X^n_{s_l})|\F^n_u] \cvL1 \E_{\P}[f(X_{s_1}, \ldots, X_{s_l})|\F_u].
\end{equation}
Using that $H$, $G$, and $f$ are continuous and bounded, we can show
that:
\begin{eqnarray}
\label{cv1}
&&\int\int f(X^{n}_{s_1}(\omega), \ldots, X^{n}_{s_l}(\omega))1_{A_2}(v)H(X^{n}_{t_1}(\omega), \ldots ,X^{n}_{t_k}(\omega)) \nonumber\\
&&      \quad \quad \quad \quad         G(\tau^{*,n}(\omega,v) \wedge u)g(v)d(\P \otimes \mu)(\omega,v)  \\
&\xrightarrow[n \to +\infty]{} &
\int\int f(X_{s_1}(\omega), \ldots, X_{s_l}(\omega))1_{A_2}(v)H(X_{t_1}(\omega), \ldots , X_{t_k}(\omega))\nonumber\\
&&      \quad \quad \quad \quad G(\tau^*(\omega,v) \wedge u)g(v)d(\P \otimes \mu)(\omega,v).
\nonumber
\end{eqnarray}

\medskip

On the other hand, $\E[f(X_{s_1}, \ldots, X_{s_l})1_{A_2}|\F_u \times \mathcal{B}]=\E[f(X_{s_1}, \ldots, X_{s_l})|\F_u]1_{A_2}$. \\
Using again that $H$, $G$ and $f$ are continuous and bounded and the
convergence (\ref{mer1}), we have:
\begin{eqnarray}
\label{cv2}
&&\int\int \E[f(X^n_{s_1}, \ldots, X^n_{s_l})1_{A_2}|\F^n_u \times \mathcal{B}](\omega, v)
        H(X^{n}_{t_1}(\omega), \ldots , X^{n}_{t_k}(\omega))\nonumber\\
&& \quad \quad  \quad \quad G(\tau^{*,n} \wedge u)g(v)d(\P \otimes \mu)(\omega, v)\nonumber \\
& \xrightarrow[n \to +\infty]{} &
\int\int \E[f(X_{s_1}, \ldots, X_{s_l})1_{A_2}|\F_u \times \mathcal{B}](\omega, v)
        H(X_{t_1}(\omega), \ldots , X_{t_k}(\omega))\nonumber \\
&& \quad \quad  \quad \quad G(\tau^*(\omega, v) \wedge u)g(v)d(\P \otimes \mu)(\omega, v).
\end{eqnarray}

\medskip

But, $H(X^{n}_{t_1}, \ldots , X^{n}_{t_k})$ is $\F^{n}_u \times 
\mathcal{B}-$measurable and \mbox{$G(\tau^{n} \wedge u)$} and both
$g(U)$, where $\forall \omega \in \Omega$, $\forall v \in [0,1]$, 
$U(\omega,v)=v$, are also $\F^{n}_u \times \mathcal{B}-$measurable, 
by continuity of $G$ and $g$. It follows that \begin{eqnarray*} 
\lefteqn{\E[\E[f(X^{n}_{s_1}, \ldots, X^{n}_{s_l})1_{A_2}|\F^{n}_u
\times \mathcal{B}]
    H(X^{n}_{t_1}, \ldots , X^{n}_{t_k})G(\tau^{n} \wedge u)g(U)] }\\
& = & \E[\E[f(X^{n}_{s_1}, \ldots, X^{n}_{s_l})1_{A_2}
    H(X^{n}_{t_1}, \ldots , X^{n}_{t_k})G(\tau^{n} \wedge u)g(U)|\F^{n}_u \times \mathcal{B}]] \\
& = & \E[f(X^{n}_{s_1}, \ldots, X^{n}_{s_l})1_{A_2}H(X^{n}_{t_1}, \ldots , X^{n}_{t_k})G(\tau^{n} \wedge u)g(U)]
\end{eqnarray*}
Identifying limits in (\ref{cv1}) and (\ref{cv2}), we obtain:
\begin{eqnarray}
\label{eqfHG}
&& \int\int \E[f(X_{s_1}, \ldots, X_{s_l})1_{A_2}|\F_u \times \mathcal{B}](\omega, v)
        H(X_{t_1}(\omega), \ldots , X_{t_k}(\omega))\nonumber \\
&& \quad \quad  \quad \quad G(\tau^*(\omega, v) \wedge u)g(v)d(\P \otimes \mu)(\omega, v) \nonumber \\
& = & \int\int f(X_{s_1}(\omega), \ldots, X_{s_l}(\omega))1_{A_2}(v)
        H(X_{t_1}(\omega), \ldots , X_{t_k}(\omega))\nonumber \\
&& \quad \quad  \quad \quad G(\tau^*(\omega, v) \wedge u)g(v)d(\P \otimes \mu)(\omega, v).
\end{eqnarray}
Then, using the approximations (\ref{alpha1}), (\ref{alpha2}),
(\ref{alpha3}), (\ref{alpha45}) and the fact that $\E[f(X_{s_1},
\ldots, X_{s_l})|\F_.]$ is a c\`adl\`ag process, we can deduce from
(\ref{eqfHG}) the equality (\ref{eqZ}):
\begin{eqnarray*}
\lefteqn{
\int\int \E[Z|\F_t \times \mathcal{B}](\omega, v)1_A(\omega)1_{\{\tau^*(\omega, v) \leqslant s\}}1_B(v)d(\P \otimes \mu)(\omega,v)}\\
& = & \int\int Z(\omega, v)1_A(\omega)1_{\{\tau^*(\omega, v) \leqslant s\}}1_B(v)d(\P \otimes \mu)(\omega, v),
\end{eqnarray*}
for every $t \in [0,T]$, for every $Z =1_{A_1 \times A_2}$, $A_1 \in \F_T$, $A_2 \in \mathcal{B}$, for every $A \in \F_t$,
for every $s \leqslant t$, for every $B \in \mathcal{B}$. \\

It follows through a monotone class argument, linearity and density 
that (\ref{eqZ}) holds whenever $Z$ is $\F_T \times 
\mathcal{B}-$measurable and integrable.

Hence, for every $t \in [0,T]$, for every $Z \in L^1(\F_T \times \mathcal{B})$,
for every $C \in \G_t \times \mathcal{B}$ (by definition of $\G_t \times \mathcal{B}$),
$$\E_{\P \otimes \mu}[\E_{\P \otimes \mu}[Z|\F_t \times \mathcal{B}]1_C]=\E_{\P \otimes \mu}[Z1_C].$$

\noindent We hence have checked (\ref{indcond}), therefore the assumption of Lemma \ref{indepcond} is filled, and we readily deduce
Lemma \ref{mingamma} from (\ref{tildestar}).
\findem

Recall now inequality (\ref{limsup-epsilon}): from the definition of
$\tau^*$ and Lemma \ref{indepcond}, whose assumption is filled as we
just have shown, it follows that
\begin{eqnarray*}
\limsup \Gamma_n(T) - \varepsilon&\leqslant &\E[\gamma(\tau^*,X_{\tau^*})]\cr &\leqslant &\tilde{\Gamma}(T)=\Gamma^*(T).
\end{eqnarray*}

As such a randomized stopping time $\tau^*$ exists for arbitrary $\varepsilon >0$, we conclude that

$$\limsup \Gamma_n(T) \leqslant \Gamma^*(T)$$
and Lemma \ref{GammaGamma*} shows now that $\Gamma^*(T)=\Gamma(T)$.
Theorem \ref{thlimsup2} is proved.
\findem

To sum up this section, under the hypothesis of Theorem
\ref{cvGamma}, we have proved the inequality $\Gamma(T) \leqslant
\liminf \Gamma_n(T)$ in Theorem \ref{thliminf}. Then, we have shown
that $\Gamma(T) \geqslant \limsup \Gamma_n(T)$ when inclusion of
filtrations $\F^n \subset \F$ (in Theorem \ref{thlimsup1}) or
convergence of filtrations $\F^n \xrightarrow{w} \F$ (in Theorem
\ref{thlimsup2}) hold, provided that Aldous' Criterion for tightness
(\ref{hypA}) is filled by the sequence $(X^n)$. At last, Theorem
\ref{cvGamma} is proved.


\section{Convergence of optimal stopping times}
\label{cvtaopt}

\begin{definition} $\tau$ is an optimal stopping time for $X$ if 
$\tau$ is a $\F-$stopping time bounded by $T$ such that 
$\E[\gamma(\tau,X_\tau)]=\Gamma(T)$. \end{definition}

Some results of existence of optimal stopping time are given for instance in \citep{S_taopt} in the case of Markov processes. \\

Now, let $(X^n)_n$ be a sequence of c\`adl\`ag processes that
converges in probability to a c\`adl\`ag process $X$. Let $(\F^n)_n$
be the natural filtrations of processes $(X^n)_n$ and $\F$ the
right-continuous filtration of $X$. We suppose again that Aldous'
Criterion for tightness (\ref{hypA}) is filled and that we have the
convergence of values in optimal stopping: $\Gamma_n(T) \to
\Gamma(T)$ (see Section \ref{intro} for the notations).

We consider, if it exists, a sequence $(\tau^n_{op})_n$ of optimal
stopping times associated to the $(X^n)$. $(\tau^n_{op})_n$ is tight
so we can find a subsequence which converges in law to a random
variable $\tau$. There are at least two problems to solve. First, is
$\tau$ a $\F-$stopping time or (at least) is the law of $\tau$ the
law of a $\F-$stopping time ? Then, if the answer is positive, is 
$\tau$ optimal for $X$, i.e. have we 
$\E[\gamma(\tau,X_\tau)]=\Gamma(T)$ ? \\

It is not difficult to answer the second question as next result 
shows: \begin{lemme} \label{opt} We suppose that $\Gamma_{n}(T) 
\xrightarrow[n \to +\infty]{} \Gamma(T)$  and that Aldous' Criterion
for tightness (\ref{hypA}) is filled. Let $(\tau^n_{op})_n$ be a 
sequence of optimal stopping times associated to $(X^n)_n$. Assume 
that $\tau$ is a stopping time such that, along some subsequence 
$\varphi$, $(X^{\varphi(n)},\tau^{\varphi(n)}_{op}) \cvl (X,\tau)$. 
Then $\tau$ is an optimal $\F-$stopping time. \end{lemme}

\demo $(X^{\varphi(n)},\tau^{\varphi(n)}_{op}) \cvl (X,\tau)$ so 
according to Proposition \ref{cvlXntar}, 
$(\tau^{\varphi(n)}_{op},X^{\varphi(n)}_{\tau^{\varphi(n)}_{op}}) 
\cvl (\tau,X_{\tau}).$ $\gamma$ is bounded and continuous, so
$\E[\gamma(\tau^{\varphi(n)}_{op},X^{\varphi(n)}_{\tau^{\varphi(n)}_{op}})] 
\xrightarrow[n \to +\infty]{} \E[\gamma(\tau,X_{\tau})].$
$(\tau^{\varphi(n)}_{op})$ is a sequence of optimal 
$(\F^{\varphi(n)})-$stopping times, so for every $n$, 
$\E[\gamma(\tau^{\varphi(n)}_{op},X^{\varphi(n)}_{\tau^{\varphi(n)}_{op}})] 
= \Gamma_{\varphi(n)}(T)$ where $\Gamma_{\varphi(n)}(T)$ is the
value in optimal stopping for $X^{\varphi(n)}$. Moreover, 
$\Gamma_{\varphi(n)}(T) \xrightarrow[n \to +\infty]{} \Gamma(T).$ 
So, by unicity of the limit, $\E[\gamma(\tau,X_{\tau})] = 
\Gamma(T).$ Finally, $\tau$ is an optimal $\F-$stopping time.
\findem

Now, it remains to find a criterion to determine wether the limit of 
a sequence $(\tau^n)$ of $(\F^n)-$stopping times is a $\F-$stopping 
time. Next proposition gives such a criterion involving convergence 
of filtrations, and which will prove useful in the applications of 
Section \ref{appl}.

\begin{prop} \label{sssuitecvf} We suppose $\F^n \xrightarrow{w} 
\F$. Let $(\tau^n)_n$ be a sequence of $(\F^n)-$stopping times that
converges in probability to a $\F_T-$measurable random variable 
$\tau$. Then $\tau$ is a $\F-$stopping time. \end{prop}

\demo $\tau^n \cvp \tau$ so $1_{\{\tau^n \leqslant .\}} \cvp 
1_{\{\tau \leqslant .\}}$ for the Skorokhod topology.\\ We fix $t$
such that $\P[\tau=t]=0$. Then, $1_{\{\tau^n \leqslant t\}} \cvp 
1_{\{\tau \leqslant t\}}.$ The sequence $(1_{\{\tau^n \leqslant
t\}})_n$ is uniformly integrable, so $1_{\{\tau^n \leqslant t\}}
\xrightarrow{L^1} 1_{\{\tau \leqslant t\}}.$ $\tau$ is
$\F_T-$measurable, so $1_{\{\tau \leqslant t\}}$ is 
$\F_T-$measurable. As $1_{\{\tau^n \leqslant t\}} \xrightarrow{L^1} 
1_{\{\tau \leqslant t\}}$, $\F^n \xrightarrow{w} \F$ and $1_{\{\tau
\leqslant t\}}$ is $\F_T-$measurable, according to \citep[Remark
2]{cvfiltration}, we have: $$\E[1_{\{\tau^n \leqslant t\}}|\F^n_.]
\cvp \E[1_{\{\tau \leqslant t\}}|\F_.].$$ \indent Let us prove that
$\E[1_{\{\tau^n \leqslant t\}}|\F^n_t] \cvp \E[1_{\{\tau \leqslant 
t\}}|\F_t].$ \\ Fix $\eta >0$ and $\varepsilon >0$. \\ $\E[1_{\{\tau
\leqslant t\}}|\F_.]$ is a c\`adl\`ag process, so we can find $s \in
]t,T]$ satisfying $\P[\Delta \E[1_{\{\tau \leqslant t\}}|\F_s] \not=
0] =0$ and such that $$\P[|\E[1_{\{\tau \leqslant t\}}|\F_s] -
\E[1_{\{\tau \leqslant t\}}|\F_t]| \geqslant \eta/3] \leqslant
\varepsilon/3.$$ Then, we have $\E[1_{\{\tau^n \leqslant
t\}}|\F^n_s] \cvp \E[1_{\{\tau \leqslant t\}}|\F_s]$ and we can find
$n_0$ such that for every $n \geqslant n_0$, $$\P[|\E[1_{\{\tau^n 
\leqslant t\}}|\F^n_s] - \E[1_{\{\tau \leqslant t\}}|\F_s]|
\geqslant \eta/3]\leqslant \varepsilon/3.$$ On the other hand,
$$\P[|\E[1_{\{\tau^n \leqslant t\}}|\F^n_t] -\E[1_{\{\tau^n 
\leqslant t\}}|\F^n_s] \geqslant \eta/3]=0$$ because $\{\tau^n
\leqslant t\} \in \F^n_t$ as $(\tau^n)_n$ is a sequence of
$(\F^n)-$stopping times, and $\{\tau^n \leqslant t\} \in \F^n_s$ 
since $s \geqslant t$. \\ Finally, for every $n \geqslant n_0$, 
\begin{eqnarray*} \lefteqn{\P[|\E[1_{\{\tau^n \leqslant t\}}|\F^n_t] 
- \E[1_{\{\tau \leqslant t\}}|\F_t]| \geqslant \eta] } \\ &
\leqslant & \P[|\E[1_{\{\tau^n \leqslant t\}}|\F^n_t]
-\E[1_{\{\tau^n \leqslant t\}}|\F^n_s] \geqslant \eta/3] +
\P[|\E[1_{\{\tau^n \leqslant t\}}|\F^n_s] - \E[1_{\{\tau \leqslant
t\}}|\F_s]| \geqslant \eta/3] \\ && \quad \quad + \P[|\E[1_{\{\tau
\leqslant t\}}|\F_s] - \E[1_{\{\tau \leqslant t\}}|\F_t]| \geqslant
\eta/3]\\ & \leqslant & \varepsilon. \end{eqnarray*} Hence,
$$\E[1_{\{\tau^n \leqslant t\}}|\F^n_t] \cvp \E[1_{\{\tau \leqslant 
t\}}|\F_t].$$ But, $(\tau^n)_n$ is a sequence of $(\F^n)-$stopping
times, so $\forall n$, $\E[1_{\{\tau^n \leqslant 
t\}}|\F^n_t]=1_{\{\tau^n \leqslant t\}}.$ Moreover, $1_{\{\tau^n
\leqslant t\}} \cvp 1_{\{\tau \leqslant t\}}.$ By unicity of the
limit, $\E[1_{\{\tau \leqslant t\}}|\F_t]= 1_{\{\tau \leqslant t\}}$ 
$a.s$. Then, for every $t$ such that $\P[\tau=t]=0$, $\{\tau 
\leqslant t\} \in \F_t.$\\

Next, the right continuity of $\F$ implies that for every $t$, 
$\{\tau \leqslant t\} \in \F_t$. Finally $\tau$ is a $\F-$stopping 
time. \findem

\begin{rem} \label{hyp_incl} A sufficient condition to get the 
$\F_T-$measurability of the limit may be the inclusion of terminal 
$\sigma-$fields $\F^n_T \subset \F_T,\forall n$. Indeed, under this 
hypothesis, $(\tau^n)$ is a sequence of $\F_T-$measurable variables. 
Hence the limit is also $\F_T-$measurable. \end{rem}

\begin{rem} Even if convergence of filtrations $\F^n \xrightarrow{w} 
\F$ holds, the limit of a sequence of $(\F^n)-$stopping times is not
a priori $\F_T-$measurable. For example, if $\F$ is the trivial 
filtration, the assumption $\F^n \xrightarrow{w} \F$ is always true. 
However, the limit of a sequence of $(\F^n)-$stopping times may not 
be a constant, so it is not always $\F_T-$measurable. \end{rem}

Proposition \ref{sssuitecvf} and Lemma \ref{opt} allow us to give a result of convergence of optimal stopping times when the processes
$X^n$ have independent increments.

\begin{theo}
\label{cv_taopt} Let $(X^n)$ be a sequence of c\`adl\`ag processes
which converges in law to a quasi-left continuous process X. We
suppose that the processes $(X^n)$ have independent increments. Let
$(\F^n)$ be the natural filtrations of processes $(X^n)$ and $\F$ be
the right-continuous filtration associated to the process $X$. Let
$(\tau^n)$ be a sequence of optimal $(\F^n)-$stopping times. If
$(X,\tau)$ is the limit in law of a subsequence of
$((X^n,\tau^n))_n$ and if $\tau$ is $\F_T-$measurable, then $\tau$
is an optimal stopping time for $X$.
\end{theo}

\demo $((X^n,\tau^n))$ is tight because $(\tau^n)$ is bounded and 
$(X^n)$ is convergent. So we can extract a subsequence 
$((X^{\varphi(n)},\tau^{\varphi(n)}))$ which converges in law to 
$(X,\tau).$ \\ Let us prove that $\tau$ is a $\F-$stopping time. \\ 
Using the Skorokhod representation theorem, we can find a 
probability space $(\tilde{\Omega},\tilde{\G},\tilde{\P})$ on which 
are defined $(\tilde{X}^{\varphi(n)},\tilde{\tau}^{\varphi(n)}) \sim 
(X^{\varphi(n)},\tau^{\varphi(n)})$ and $(\tilde{X},\tilde{\tau})
\sim (X,\tau)$ such that
$(\tilde{X}^{\varphi(n)},\tilde{\tau}^{\varphi(n)}) 
\xrightarrow{a.s} (\tilde{X},\tilde{\tau})$ in
$(\tilde{\Omega},\tilde{\G},\tilde{\P}).$ As $(X^n)$ are processes 
with independent increments, $(\tilde{X}^{\varphi(n)})$ also are. 
Using \citep[Proposition 3]{Memin_2003}, we have the extended 
convergence $$(\tilde{X}^{\varphi(n)},\F^{\tilde{X}^{\varphi(n)}}) 
\cvp (\tilde{X},\F^{\tilde{X}})$$ where
$\F^{\tilde{X}^{\varphi(n)}}$ (resp. $\F^{\tilde{X}}$) is the 
natural filtration of the process $\tilde{X}^{\varphi(n)}$ (resp. 
$\tilde{X}$). \\ On the other hand, $\tau$ is $\F_T-$measurable. So, 
we can find a measurable function $f$ such that $\tau=f(X)$. 
$(\tilde{X},\tilde{\tau}) \sim (X,\tau)$ and $(X,\tau)=(X,f(X))$, so 
$\tilde{\tau}=f(\tilde{X})$ $a.s$. Hence, $\tilde{\tau}$ is 
$\F^{\tilde{X}}_T-$measurable. \\ Moreover, 
$\tilde{\tau}^{\varphi(n)} \xrightarrow{a.s.} \tilde{\tau}$ by 
construction and, as $(\tau^{\varphi(n)})$ is a sequence of 
$(\F^{\varphi(n)})-$stopping times, $(\tilde{\tau}^{\varphi(n)})$ is 
a sequence of $(\F^{\tilde{X}^{\varphi(n)}})-$stopping times. \\ 
Then, using Proposition \ref{sssuitecvf}, $\tilde{\tau}$ is a 
$\F^{\tilde{X}}-$stopping time. \\ Next, using \citep[Proposition 
16.20]{preprintAldous}, Aldous' Criterion for tightness is filled
because $\tilde{X}$ is quasi-left continuous and 
$(\tilde{X}^{\varphi(n)},\F^{\tilde{X}^{\varphi(n)}}) \cvl 
(\tilde{X},\F^{\tilde{X}})$. Moreover,
$\Gamma^{\tilde{X}^{\varphi(n)}}(T) \to \Gamma^{\tilde{X}}(T)$ 
according to Theorem \ref{cvGamma}. \\ Then, according to Lemma 
\ref{opt}, $\tilde{\tau}$ is an optimal $\F^{\tilde{X}}-$stopping 
time.\\ Finally, $\tau$ is an optimal $\F-$stopping time. \findem


\section{Applications}
\label{appl}

\subsection{Application to discretizations}

\begin{prop}
\label{cvGammadiscr}
Let us consider a quasi-left continuous process X with independent increments.
Let \mbox{$(\pi^n=\{t_1^n, \ldots t_{k^n}^n\})_n$} be an increasing sequence of subdivisions of $[0,T]$ with mesh going to 0
$(|\pi^n| \xrightarrow[n \to +\infty]{} 0)$. We define the sequence of discretized processes $(X^n)_n$ by $\forall n$, $\forall t$,
\mbox{$X_t^n = \sum_{i=1}^{k^n-1} X_{t_i^n} 1_{t_i^n \leqslant t < t_{i+1}^n}$.}  \\
Let us denote by $\F$ the right-continuous natural filtration of $X$ and by $(\F^n)_n$ the natural filtrations of the $(X^n)_n$. \\
Then, using the notations of Section \ref{intro}, \mbox{$\Gamma_n(T) \xrightarrow[n \to +\infty]{} \Gamma(T)$.}
Moreover, if $(\tau^n)$ is a sequence of stopping times associated to the processes $(X^n)$ and if $(X,\tau)$ is the limit in law of a
subsequence of $(X^n,\tau^n)$, then $\tau$ is an optimal stopping time for $X$.
\end{prop}

\demo $X^n \xrightarrow[n \to +\infty]{} X$ $a.s.$ then in 
probability, for every $n$ $\F^n \subset \F$ by definition of $X^n$. 
Moreover, $X$ is quasi-left continuous and $(X^n)$ is a sequence of 
discretized processes, so we can easily check that Aldous' Criterion 
is filled. So, using Theorem \ref{cvGamma}, $\Gamma_n(T) 
\xrightarrow[n \to +\infty]{} \Gamma(T)$. \\ On the other hand, for
every $n$, $\F^n_T \subset \F_T$. So $\tau$ is $\F_T-$measurable. 
Then, according to Theorem \ref{cv_taopt}, $\tau$ is an optimal 
stopping time for $X$. \findem

\subsection{Application to financial models}

\subsubsection{The models}

The convergence of properly normalized Cox-Ross-Rubinstein models to a Black-Scholes model is a standard in financial mathematics. By
convergence, it is usually meant here convergence of option prices. We are going to apply our results to prove that a sequence of so-called
rational times of exercise for an american put in a Cox-Ross-Rubinstein converge, under the same normalization, to a rational time of
exercise for an american put in the Black-Scholes model.

We just recall here the classical notation for both models.

The Black-Scholes model on an interval $[0,T]$ consists in a market
with one non-risky asset of price $S^0_t=S^0_0e^{rt}$ at time $t$,
$r$ denoting the instant interest rate, and a risky asset whose
price is governed by the following stochastic differential equation:
\begin{equation}
S_t=S_t(\mu dt + \sigma dB_t)
\end{equation}
 where $\mu$ and $\sigma$ are positive reals and $(B_t)$ is a standart Brownian motion. We denote by $\P^*$ the risk-neutral probability,
 under which the actualized price of the risky asset is a martingale, and by $(\F_t)_{t\leq T}$ the filtration generated by the Brownian
 motion $B$.

If we are given an american put option with maturity $T$ and strike price $K$, then its optimal value is defined as
$$\Gamma^{S}(T)=\underset{\tau \in \mathcal{T}_T}{\sup} \E_{\P^*}[e^{-r\tau}(K-S_{\tau})^+],$$
where $\mathcal{T}_T$ is the set of $\F-$stopping times bounded by $T$, and the expectation is taken under $\P^*$. A rational exercise time
is then a stopping time $\tau^0$ such that
$$ \E_{\P^*}[e^{-r\tau^0}(K-S_{\tau^0})^+]=\Gamma^{S}(T).$$ \\

We now build a sequence of random walks approaching $B$, following
the construction of \citep{Knight_62}. We refer to \citep{ItoMcKean}
for explicit details. We only need to know here that Knight has
built an array $(Y^n_i)$ such that, for every $n$, $(Y^n_i)_i$ is a
sequence of $\F_T-$mesurable independent Bernoulli variables, such
that $\P[Y^n_i=1]=\P[Y^n_i=-1]=1/2$, and for which, if we put $B^n_t
=\sqrt{\frac{T}{n}} \sum_{i=1}^{[\frac{nt}{T}]} Y^n_i$, holds the
following convergence: \begin{equation} \label{cv_knight} \P \left[ 
\lim_{n \uparrow +\infty} \sup_{t \in [0,T]} |B^n_t -
B_t|=0\right]=1. \end{equation}

The last step is to build the Cox-Ross-Rubinstein models based upon 
the array $(Y^n_i)$ in such a way that holds the convergence of 
binomial prices for the risky assets (denoted by $S_n$) to $S$ (the 
reader will find the appropriate normalizations, e.g. in \citep{LL} 
or \citep{Shi_finance} or any textbook on mathematical finance). \\ 
For each $n$, the maximal expectation of profit for the associated 
Cox-Ross-Rubinstein model is given by: $$\Gamma^{S^n}(T) = 
\sup_{\tau^n \in \mathcal{T}^n_T} \E_{\P^{*,n}}[(1+rT/n)^{-([\tau^n
n/T])}(S^n_{\tau^n}-K)^+] $$ where $(\F^n)$ denotes the (piecewise
constant) filtration generated by the price process $(S^n_i)_i$ 
(which is also the filtration generated by the process $B^n$), 
$\mathcal{T}^n_T$ is the set of $\F^n-$stopping times bounded by 
$T$, and $\P^{*,n}$ is the equivalent probability making the 
actualised price process a $\F^n-$martingale.

\subsubsection{Convergence of values in optimal stopping}

Having used Knight's construction ensures us that
$$(B^n,S^n) \cvps (B,S)$$
and as the $B^n$'s are processes with independent increments and
$S^n$ and $S$ are bijective functions of $B^n$ and $B$,
\citep[Proposition 3]{Memin_2003} gives the extended convergence:
$(S^n,\F^{S^n}) \cvp (S,\F^S)$. Thanks to Theorem \ref{cvGamma^n},
whose hypothesis is clearly fulfilled, we deduce then
\begin{equation}\label{cvgamma}
\Gamma^{S^n}(T) \xrightarrow[n \to +\infty]{} \Gamma^{S}(T)
\end{equation}
\indent
Remark at last that (\ref{cvgamma}) holds regardless of the specific construction of the prelimit Cox-Ross-Rubinstein models. Indeed,
according to Remark \ref{GammaX}, the value in optimal stopping only depends on the law of the underlying process hence every
Cox-Ross-Rubinstein model with the same law as $S$ (and any Black-Scholes limiting model) would perfectly fit, provided that the
correct normalizations are performed in order to have the convergence (in law) of the price processes.\\
\indent To sum up, we have just proved the following result:
\begin{prop}
\label{cv_reduites_BS_CRR}
When approximating a Black-Scholes model by a sequence of Cox-Ross-Rubinstein  models, we have convergence of the associated sequence of
values in optimal stopping:
$$\Gamma^{S^n}(T) \xrightarrow[n \to +\infty]{} \Gamma^{S}(T).$$
\end{prop}

This convergence is already well known (see \citep{MP}, \citep{Lamberton93} or \citep{AK94} for example). What is new in this paper is the
result of convergence of optimal stopping times proved in the next section.

\subsubsection{Convergence of optimal stopping times}

In the same framework as above, let us end with the study of the 
convergence of optimal stopping times. \\ \indent $S^n$ is a Markov 
process so there exists a sequence $(\tau^n_{op})$ of optimal
$\F^n-$stopping times. The sequence $(B,B^n,S^n, \tau^n_{op})$ is 
tight, and up to some subsequence, we can assume that $$(B,B^n,S^n, 
\tau^n_{op}) \cvl (B,B,S,\tau)$$ for some random variable $\tau$.\\
By Skorokhod's representation lemma, we can assume that this 
convergence holds almost surely (while preserving the links between 
the Brownian motion $B$ and all other processes under 
consideration).

As in previous subsection, convergence of filtrations holds, 
moreover our specific construction (Knight's one) ensures that for 
every $n$, $\tau^n$ is $B^n_T-$measurable, hence $\F_T-$measurable. 
>From Proposition \ref{sssuitecvf}, we deduce that $\tau$ is a 
$\F-$stopping time.\\ As moreover, $\Gamma^{S^n}(T) \to 
\Gamma^{S}(T)$, Lemma \ref{opt} now says that $\tau$ is indeed an
optimal $\F-$stopping time. \\

We have just proved the following result:
\begin{prop}
\label{cv_taopt_BS_CRR}
When approximating a Black-Scholes model by a sequence of Cox-Ross-Rubinstein  models based on Knight's construction, if a subsequence of
$((B^n,S^n, \tau^n_{op}))_n$ -where $(\tau^n_{op})_n$ is a sequence of optimal stopping times for the prelimit models- converges in law to
$(B,S,\tau)$, then $\tau$ is an optimal stopping time for Black-Scholes model.
\end{prop}

\begin{rem}
We stress once more on the fact that, whereas Proposition \ref{cv_reduites_BS_CRR} remains true for every Cox-Ross-Rubinstein approximation
of a Black-Scholes models, the proof of Proposition \ref{cv_taopt_BS_CRR} rely upon the fact that $\tau$ is actually a stopping time for the
natural filtration associated filtration of Black and Scholes model (and not for a larger one like in the existing papers), for which we need
Knight's construction of the prelimit models.
\end{rem}



\end{document}